\documentclass[11pt,onecolumn]{article}
\usepackage{amscd}
\usepackage{times}
\usepackage{amsmath}
\usepackage{amssymb}
\usepackage{xspace}
\usepackage{theorem}
\usepackage{graphicx}
\usepackage{ifpdf}
\usepackage{url,hyperref}
\usepackage{latexsym}
\usepackage{euscript}
\usepackage{xspace}
\usepackage[all]{xy}
\usepackage{color}
\usepackage{makeidx}
\usepackage{tikz}
\long\def\remove#1{}
\newtheorem{theorem}{Theorem}[section] 

\newtheorem{obs}[theorem]{Observation}
\newtheorem{corollary}[theorem]{Corollary}

\newtheorem{definition}[theorem]{Definition}
\newtheorem{proposition}[theorem]{Proposition}

\newcommand {\mm}[1] {\ifmmode{#1}\else{\mbox{\(#1\)}}\fi}

\newcommand{\img}{\mathrm img}
\newcommand{\supp}{\mathrm supp}

\newcommand{\cancel}[1]

\begin{document}

\title{New invariants for a real valued and angle valued map.
\newline (an Alternative to Morse- Novikov theory) }

\author{
Dan Burghelea  \thanks{
Department of Mathematics,
The Ohio State University, Columbus, OH 43210,USA.
Email: {\tt burghele@math.ohio-state.edu}}
}
\date{}

\date{}
\maketitle

\begin{abstract}
\vskip .2in

This paper  but section 6  is essentially my lecture at  The Eighth Congress of Romanian Mathematicians, 
 2015, Iasi, Romania.
The paper summarizes  the definitions and the properties of the invariants  associated to a real or an angle valued map in the framework of what we call 
an Alternative to Morse--Novikov theory. These invariants  are configurations of points in the complex plane, 
configurations of vector spaces or modules  
indexed by complex numbers  and collections of Jordan cells. 
 The first are refinements of Betti numbers, the second of homology and the third of monodromy.
 Although not discussed in this paper but discussed in works this report is based on, these invariants are computer friendly (i.e. can be calculated by computer implementable algorithms when the source of the map is a simplicial complex and the map is simplicial) and are of relevance for the  dynamics of flows which admit Lyapunov real or angle valued map.
\end{abstract}
\thispagestyle{empty}
\setcounter{page}{1}

\tableofcontents

\section{Introduction}

This paper  but section 6  is essentially my lecture at  The Eighth Congress of Romanian Mathematicians, 
 2015, Iasi, Romania.

Classical Morse theory and Morse--Novikov theory consider a Riemannian manifold $(M,g)$ and a Morse real valued or a Morse angle valued map, $f:M\to \mathbb R$ or $f:M\to \mathbb S^1,$  and relate the  dynamical invariants of the vector field $grad_g \ f,$  namely

-- the  rest points of $grad_g f$= critical points of $f,$

-- the  instantons \footnote {isolated trajectories between critical points} between two rest points  $x,y$ of $grad_g f,$

-- the  closed trajectories of $grad_g f$ (when $f$ is angle valued) 

\noindent to the algebraic topology of the underlying manifold $M$ or of the pair $(M, \xi_f)$ in case $f$ is angle valued map. Here $\xi_f$ denotes the  degree one integral cohomology class represented by $f.$

The results of the theory can be applied to any vector field $V$ on $M$  which admits a closed differential one form $\omega \in \Omega^1(M)$ as Lyapunov rather than $grad_g f$, since  the dynamics of such vector field $V$ (when generic) is the same as of $grad_g f$ for some Riemannian metric $g$ and some $f,$ angle valued map cf \cite{BH08}. 
The results of the theory can be used in both ways; knowledge of the dynamical invariants of $grad_g f$ permits to calculate the topological invariants of $M$ or of $(M,\xi_f)$ and the algebraic topological invariants of $M$ or of $(M,\xi)$ provide significant constraints for dynamics of a vector field with Lyapunov map representing $\xi,$ cf \cite {BH08}.

\vskip .1in 

The ANM theory associates to a pair $(X,f)$, $X$ a compact ANR,  $f$ a continuous real or angle valued map defined on $X$ and  $\kappa$ a field  a collection of invariants: the configurations $\delta^f_r, \ \hat \delta^f_r, \ \hat{\hat \delta}^f_r$ and the Jordan cells $\mathcal J_r(f), r\geq 0.$  

 The configuration $\delta^f_r$ is a finite collection of  points with multiplicity located in $\mathbb C$ in case $f$ is  real valued  and in $\mathbb C\setminus 0$ in case $f$ is   angle valued and the configuration $\hat \delta^f_r$ is given by the same points but instead of natural numbers as multiplicities have $\kappa-$vector spaces or free $\kappa[t^{-1},t]-$modules  assigned to, where $\kappa[t^{-1},t]$ denotes the ring of Laurent polynomials with coefficients in $\kappa.$
 A Jordan cell is a pair $(\lambda, k)$  with $\lambda$ a nonzero element in the algebraic closure of the field $\kappa$ and $k$  a positive integer. The pair $(\lambda,k)$ is an abbreviation  for the $k\times k$ Jordan matrix 
$$T(\lambda, k)= \begin{pmatrix} \lambda&1&0&0&\cdots&0\\
0&\lambda&1 &0&\cdots &0\\ 
\cdots \\
0&0&\cdots& 0&\lambda& 1\\
0&0&\cdots &0&0&\lambda\end{pmatrix}.$$

The configurations $\delta^f_r$ and the collections $\mathcal J_r(f), r\geq 0$ are {\it computer friendly} in the sense that for a simplicial complex and a simplicial map can be calculated by computer implementable algorithms.

On one side these invariants refine basic algebraic topology invariants of $X$ and $(X;\xi_f)$ (Betti numbers or Novikov-Betti numbers, Homology or Novikov  homology, monodromy).
On other side they are closed to the dynamical elements (rest points, instantons, closed trajectories) of a flow on $X$ which has $f$ as a Lyapunov map and permit to detect the presence and get informations about the cardinality of such elements.

The configuration $\delta^f_r$ is a configuration of points in the complex plane, each such point 
corresponding to a pair of critical values of $f$ (i.e. bar codes in the terminology of \cite {BD11}) whose multiplicity have homological interpretation. The configuration $\hat \delta^f$ is a configuration of vector spaces or modules indexed by complex numbers with the vector space or module  $\hat \delta^f_r(z)$ of dimension or rank
equal to $\delta^f_r(z)$ and specifying a piece of the  homology $H_r(X)$ or Novikov homology $H^N_r(X, \xi_f).$  The Jordan cells 
$\mathcal J_r(f)$ are pairs $(\lambda, k)$ each providing a Jordan matrix which appears in the Jordan decomposition of the $r-$ monodromy of $\xi_f.$

In contrast with the classical Morse--Novikov theory concerned with critical points of $f,$  instantons  and periodic orbits of $grad_g f$ for $X$ a smooth manifold and $f$ a Morse real  or angle valued map, the configurations $\delta^f_r,$ $\hat \delta^f_r$ and the Jordan cells $\mathcal J_r,$  associated to $f$  in  AMN-theory, 

--\   are defined for  spaces $X$ and maps $f$ considerably more general than manifolds and Morse maps,

--\  are computable by effective algorithms  when $X$ is a finite simplicial complex and $f$ simplicial map,
  
--\  enjoy robustness to $C^0-$ perturbation and satisfy Poincar\'e duality. 

This paper summarizes the definitions and the properties of  the  invariants $\delta^f_r, \hat\delta^f_r, \hat{\hat \delta}^f_r, \mathcal J_r(f) $ in AMN-theory  and addresses only the first aspect of the theory, the algebraic topology aspect.  It also indicates a few mathematical  applications (section 6). 
The results are stated in Sextion 4. Details for the proofs are contained in \cite {B1}, \cite {B2}, \cite {B3} and partially in \cite {BH} where the computational aspects  of these invariants are also addressed.

\section{Preliminary definitions}

\subsection {Configurations}\label {C}

Let $X$ be a topological space and $\kappa$ a fixed field. 
A configuration of points in $X$ is a map $\delta:X\to \mathbb Z_{\geq 0}$ with finite support
and  a configuration of $\kappa-$vector spaces or of free $\kappa[t^{-1},t]-$modules indexed by the points in $X$ is a map $\hat \delta$ defined on $X$ with values $\kappa-$vector spaces or free $\kappa[t^{-1},t]-$modules with finite support.  A point $x\in X$ is in the support of $\delta$ if $\delta (x)\ne 0$ and in the support of $\hat\delta$ if $\hat\delta(x)$ is of dimension or of rank different from $0.$ The non negative integer $\sum_{x\in X}  \delta(x)$ is referred to as the {\it cardinality} of $\delta.$ One denotes by $\mathcal C_N(X)$ the set of configurations of cardinality  $N.$ 

One says that the configuration $\hat\delta$ refines  the configuration $\delta$ if $ \dim \hat\delta(x)= \delta(x).$   

If $\kappa= \mathbb C$ one can consider also  configurations with values in Hilbert modules of finite type over a von Neumann algebra, in our discussion  always $\mathbb L^\infty(\mathbb S^1),$ the finite von Neumann algebra obtained by the von Neumann completion of the group ring $\mathbb C[\mathbb Z]$ which is exactly $\mathbb C[t^{-1},t].$

Let $V$ be a finite dimensional vector space over $\kappa$ a field or a free f.g. $\kappa[t^{-1},t]-$ module or  a finite type Hilbert module over $L^\infty (\mathbb S^1).$ 
Consider the set $\mathcal P(V)$  of  subspaces of $V,$ split free submodules of $V,$ closed Hilbert submodules of $V$ respectively. 
 
One denotes by $\mathcal C_V(X) $ the set of configurations with values in $\mathcal P(V)$ which satisfy the property that 
the induced map $I_\delta: \oplus _{x\in X} \hat\delta(x)\to V$ is an isomorphism. An element of 
$\mathcal C_V(X) $ will be denoted by $\hat{\hat \delta}$ rather than $\hat \delta$ to emphasize the additional properties.

The sets $\mathcal C_N(X)$ and ${\bf \mathcal C}_V (X)$ carry  natural topologies, referred to as the {\it collision topology}.
One  way to describe these topologies is to specify for each $\delta$ or $\hat\delta$  a system of {\it fundamental neighborhoods}.

If $\delta$ has as support  the set of points $\{x_1, x_2, \cdots  x_k\},$ a fundamental neighborhood $\mathcal U$ of $\delta$ is specified by a collection of $k$ disjoint open neighborhoods  $U_1, U_2,\cdots, U_k$ of $x_1,\cdots  x_k,$ and consists of $\{\delta'\in \mathcal C_N(X)\mid
\sum_{x\in U_i} \delta'(x)=\delta(x_i)\}.$ 
Similarly if $\hat{\hat \delta}$ has as support  the set of points $\{x_1, x_2, \cdots  x_k\}$ with $\hat{\hat\delta}(x_i)= V_i\subseteq V\},$ a fundamental neighborhood $\mathcal U$  of $\hat{\hat \delta}$ is specified by a collection of $k$ disjoint open neighborhoods  $U_1, U_2,\cdots U_k$ of $x_1,\cdots  x_k,$ and consists of configuration $\hat{\hat \delta}'$ which satisfy the following: 

a)  for any $x\in U_i$ one has $\hat{\hat\delta}'(x)\subset V_i,$

b)  the map $I_{\hat{\hat\delta}'} ( \oplus_{x\in U_i} \hat {\hat \delta}' (x))= V_i.$ 

\noindent 
Note that 
\begin{obs}\label {O211}\ 

\begin{enumerate}
\item $\mathcal C_N(X)$ identifies to  the  $N-$fold symmetric  product $X^N/\Sigma_N$  of $X$ \footnote { $\Sigma _N$ is the group of permutations of $N$ elements} and if $X$ is  a metric space with distance $D$ then  
the collision topology  is the same  as the topology defined by  the metric $\underline D$ on $X^N/\Sigma_N$  induced from the distance $D.$   This induced metric is referred to as the {\it canonical metric} on  $\mathcal C_N(X).$
\item If $X= \mathbb C$ then $\mathcal C_N(X)$ identifies to the degree $N-$monic polynomials with complex coefficients and if $X= \mathbb C\setminus 0$  to the degree $N-$monic polynomials  with non zero free coefficient.    
To the configuration $\delta$ whose support consists of the points $z_1, z_2, \cdots z_k$ with  $\delta(z_i)= n_i$ one associates the monic polynomial  $P^\delta (z)= \prod _i (z- z_i)^{n_i}. $ 
Then as topological spaces $\mathcal C_N(\mathbb C)$ identifies to $\mathbb C^N$  and $\mathcal C_N(\mathbb C\setminus 0)$ to $\mathbb C^{N-1}\times (\mathbb C\setminus 0).$ 
\item 
If $X= \mathbb T:=\mathbb R^2/\mathbb Z,$ the quotient of $\mathbb R^2$by the action $ \mu (n, (a,b))= (a+2\pi n, b+2\pi n),$  the space $\mathbb T$ can be identified to $\mathbb C\setminus 0$ by $\langle a, b\rangle \rightarrow e^{ia+(b-a)}$ then 
$\mathcal C_N(\mathbb T)$ and $\mathcal C_N(\mathbb C\setminus 0)$ are homeomorphic. Here  $\langle a, b\rangle$ denotes the $\mu-$orbit of $(a,b).$  
\item The {\it canonical metrics} $\underline D$ on $\mathcal C_N(\mathbb R^2)$ or 
$\mathcal C_N(\mathbb T)$referrs to the metric derived from the complete Euclidean metric $D$ on $\mathbb R^2$ or $\mathbb R^2/ \mathbb Z.$  Both these  metrics are complete. 
Note that standard metric on $\mathbb C\setminus 0$ is not complete so although $\mathbb T$ and $\mathbb C\setminus 0$ are homeomorphic,  hence so are $\mathcal C_N(\mathbb T)$ and 
 $\mathcal C_N(\mathbb  C\setminus 0),$  when equipped with the canonical metric they are not isometric.  
\end{enumerate}
\end{obs}

\subsection {Tame maps}\label {SS22}

A space $X$ is an ANR if any closed subset $A$ of a metrizable space $B$ homeomorphic to $X$ has a neighborhood $U$ which retracts to $A$, cf \cite {Hu} chapter 3. 
Any space homeomorphic to a locally finite simplicial complex  or to a finite dimensional topological manifold or an infinite dimensional manifold (i.e. a paracompact separable Hausdorff space locally homeomorphic to the infinite dimensional separable Hilbert space  or to the Hilbert cube $[0,1]^\infty$\ \footnote{product of county;e copies of the interval $[0,1]$}) is an ANR. 

\begin{enumerate}
\item A continuous proper map $f:X\to \mathbb R,$ $X$ an ANR  \footnote {This rules out infinite dimensional Hilbert manifolds} is {\it weakly tame}  if for any $t\in \mathbb R,$  the level $f^{-1}(t)$ is an ANR. Therefore for any bounded or unbounded closed interval $I$ the space $f^{-1}(I)$ is an ANR. 
\item The number $t\in \mathbb R$ is a {\it regular value} 
if there exists $\epsilon >0$ small s.t. for any $t'\in (t-\epsilon, t+\epsilon)$  the inclusion $f^{-1}(t' )\subset f^{-1}(t-\epsilon, t+\epsilon)$ is a homotopy equivalence. A number $t$ which is not regular value is a {\it  critical value}. 
In different words the homotopy type  of the $t-$level does not change in the neighborhood of a regular value  and does change in any  neighborhood of a critical value. One denotes by $Cr(f)$ the collection of critical values of $f.$
 \item The map $f$ is called {\it tame} if weakly tame and in addition: 

i) The set of  critical values $Cr(f)\subset \mathbb R$ is discrete, 
 
ii) The number  $\epsilon (f):= \inf \{|c-c'| \mid c,c'\in Cr(f), c\ne c'\}$ satisfies $\epsilon(f)>0.$  

If $X$ is compact then (i)  implies (ii).

\item An ANR for which the set of 
tame maps is dense in the space of all maps w.r. to the fine- $C^0$ topology is called a {\it good ANR}.

There exist compact ANR's (actually compact homological n-manifolds) 
with no co-dimension one subsets which are ANR's, hence compact ANR's which are 
not {\it good } , cf \cite {DW}. 
\end{enumerate}
The reader should be aware of the following rather obvious facts.
\begin{obs}\label {O21}\
\begin{enumerate}
\item If $f$ is a weakly tame  map then the compact ANR $f^{-1}([a,b])$ has the homotopy type of a finite simplicial complex (cf \cite{Mi2}) and therefore has finite dimensional homology w.r. to any field $\kappa.$
\item If $X$ is a locally finite simplicial complex and $f$ is linear on each simplex then $f$ is weakly tame  with the set of critical values discrete. Critical values are among the values $f$ takes on vertices. If in addition $X$ is compact then $f$ is tame.
If $M$ is a smooth manifold and $f$ is proper smooth map with all critical points of finite codimension, in particular $f$ is a Morse map, then $f$ is weakly tame and when  $M$ is compact $f$ is tame.  
\item If $X$ is homeomorphic to a compact simplicial complex or to a compact topological manifold the set of tame maps 
is dense in the set of all continuous maps equipped with the $C^0-$topology (= compact open topology). The same remains true if $X$ is a compact Hilbert cube manifold defined in the next section. In particular all these spaces are good ANR's.
\item On a smooth manifold  the Morse functions are dense in the space of all continuous function 
w.r. to the fine $C^0-$topology and are generic in any $C^r-$ topology, $r\geq 2.$
\end{enumerate}
\end{obs}

\subsection{Algebraic topology}

Let $\kappa$ be a field.
For an ANR $X$ denote by $H_r(X)$ the (singular) homology with coefficients in $\kappa;$ this is a $\kappa-$vector space which when $X$ is compact is finite dimensional by \cite{Mi2}.

Denote by $\beta_r(X):=\beta_r(X;\kappa)= \dim H_r(X)\ r \geq 0$  referred below as the $r-$th Betti number and by $\chi(X)= \chi(X;\kappa)=\sum_r (-1)^r \beta_r(X)$ the Euler characteristic with coefficients in $\kappa.$

For a pair $(X, \xi\in H^1(X;\mathbb Z)),$  $X$ a compact ANR and $\xi$ a degree one  integral cohomology class,  consider $\pi: \tilde X\to X$ an infinite cyclic cover associated to $\xi$ (unique up to isomorphism), and let $\tau:\tilde X\to \tilde X$ be the generator of the group of deck transformations (the infinite cyclic group $\mathbb Z ).$ 

The space $\tilde X$ is a locally compact ANR and the $\kappa-$vector space $H_r(\tilde X)$ is a finitely generated $\kappa[t^{-1},t]-$ module with the multiplication by $t$ given by the isomorphism $T_r:H  (\tilde X)\to H_r(\tilde X)$ induced by the homeomorphism $\tau.$ The submodule of torsion elements of $H_r(\tilde X),$  denoted by $V_r(X;\xi),$ when regarded as a $\kappa-$vector space is  finite dimensional and the $\kappa[t^{-1},t]-$module $ H_r(\tilde X)/ V_r(X;\xi)$ is free 
of finite rank.

The isomorphism class of the $\kappa[t^{-1},t]-$module $V_r(X;\xi),$  equivalently of the pair $(V_r(X;\xi), T_r)$ with $V_r(X;\xi)$ viewed as a $\kappa-$vector space with a linear automorphism  $T_r,$  is referred to as the $r-$th monodromy. The free $\kappa[t^{-1},t]-$module $H^N_r(X:\xi):= H_r(\tilde X)/ V_r(X;\xi)$  is referred below as the Novikov homology in dimension $r,$ and  its rank as the 
$r-$Novikov--Betti number and denoted by $\beta^N_r(X;\xi).$

If $\kappa =\mathbb C$ is the field of complex numbers then the ring $\mathbb C[t^{-1},t],$ equivalently the group algebra $\mathbb C[\mathbb Z],$  has a canonical completion to the finite von-Neumann algebra $L^\infty(\mathbb S^1)$ and the module $H^N_r( X;\xi)$ to a finite type $L^\infty(\mathbb S^1)-$Hilbert module,  of  von-Neumann dimension  $\beta^B_r(X;\xi).$  The completion of $H^N_r( X;\xi)$is exactly the $L_2-$homology 
$H^{L_2}_r(\tilde  X),$ 
cf \cite{B2}.  The completion  of $H^N_r(X;\xi)$ is referred to as the von Neumann completion, 
and depends a priory on additional data 
such as:  a Riemannian metric when  $X$ a compact smooth manifold, a triangulation when $X$ is a finite simplicial complex or more algebraically, an inner $\mathbb C[t^{-1}, t]-$product on $H^N_r(X;\xi)$, but all these data  lead to isomorphic $L^\infty(\mathbb S^1)-$Hilbert modules, 
cf\cite{B2}.
\vskip .1in 
\section {The configurations 
and the set of Jordan cells 
}

Let $f:X\to \mathbb R$ be a  a proper continuous map, $X$ an ANR and  $\kappa$ be a fixed field. Denote by: 

\vskip .1in 
--\ $X_a,$ the sub level 
$X_a: =f^{-1}(-\infty,a]),$ 

--\ $X^b,$ the super level 
$X^b:=f^{-1}([b,\infty)),$
 
--\ $\mathbb I^f_a(r):= \img (H_r(X_a)\to H_r(X)) \subseteq H_r(X),$

--\ $\mathbb I^b_f(r):= \img (H_r(X^b)\to H_r(X))\subseteq H_r(X),$

--\ $\mathbb F_r^f(a,b):= \mathbb I^f_a(r)\cap \mathbb I^b_f(r)\subseteq H_r(X),$ $F^f_r(a,b)= \dim \mathbb F^f_r(a,b).$

\begin{obs}\label {O31}\

1.\ If $a'\leq a$ and $b\leq b'$ then $\mathbb F_r^f(a',b')\subseteq\mathbb F_r^f(a,b).$

2.\ If $a'\leq a$ and  $b\leq b'$ then $\mathbb F_r^f(a',b)\cap \mathbb F_r^f(a,b')= \mathbb F_r^f(a',b').$

3.\ $\dim F_r(a,b)<\infty$ (cf \cite {B1} Proposition 3.4).

4.\ $\sup_{x\in X} |f(x)-g(x)|
<\epsilon$ implies $\mathbb F^g(a-\epsilon, b+\epsilon) \subseteq\mathbb F^f_r(a,b).$

5.\ If $f$ is weakly tame and number $a\in \mathbb R$ is a {\it regular value} then there exists $\epsilon>0$ so that for any $0\leq t, t'<\epsilon$
the inclusions $\mathbb I^f_{(a-t)}(r) \subseteq \mathbb I^f_{(a+t')}(r)$ and  $\mathbb I_f^{(a-t')}(r) \supseteq \mathbb I_f^{(a+t)}(r)$ are isomorphisms for all $r$.
\end{obs}
\vskip .1in

A set $B\subset \mathbb R^2$  of  the form $B= (a',a]\times [b,b')$  with $a'<a, b<b'$ is called  {\it  box}. For $(a,b)\in \mathbb R^2$ and $\epsilon >0$ denote by $B(a,b;\epsilon) $ the box $B(a,b;\epsilon): = (a-\epsilon, a]\times [b, b+\epsilon).$
 To the box $B$  we assign the vector space $$\mathbb F^f_r(B):= \mathbb F^f_r(a,b)/ \mathbb F^f_r(a',b)+\mathbb F^f_r(a,b')$$ of dimension 
 $$F^f_r(B): = \dim \mathbb F^f_r(B).$$ 
In view  of Observation \ref{O31} item 3 \ $F^f_r(B) <\infty $ and in view of  
Observation \ref{O31} item 2  

\hskip 1in $F^f_r(B):=  F^f_r(a,b) + F^f(a',b') - F^f_r(a',b)-  F^f(a,b') .$

\vskip .1in 

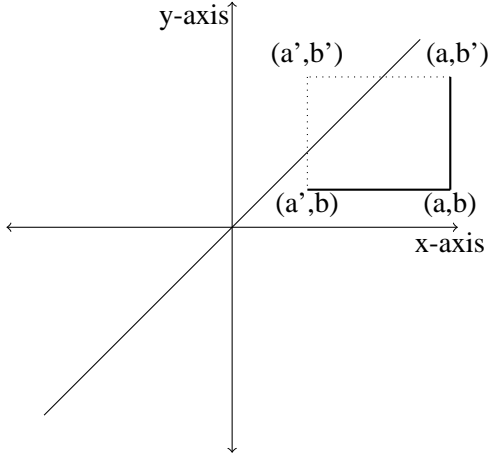
\begin{figure}
\begin{tikzpicture}
\draw [<->]  (0,3) -- (0,0) -- (3,0);
\node at (-0.5,2.8) {y-axis};
\node at (2.9,-0.2) {x-axis};
\node at (1,0.3) {(a',b)};
\node at (2.9,0.3) {(a,b)};
\node at (3,2.3) {(a,b')};
\node at (1,2.3) {(a',b')};
\draw [<->]  (0,-3) -- (0,0) -- (-3,0);
\draw [thick] (1,0.5) -- (2.9,0.5);
\draw [thick] (2.9,0.5) -- (2.9,2);
\draw [dotted] (2.9,2) -- (1,2);
\draw [dotted] (1,2) -- (1,0.5);
\draw (0,0) -- (2.5,2.5);
\draw (0,0) -- (-2.5,-2.5);
\end{tikzpicture}
\caption {The {\it box} $ B : =(a',a]\times [b,b')\subset \mathbb R^2$  }
\end{figure}

\newcommand{\mynewnewpicture}[1][ ]{
\begin{tikzpicture} [scale=1]
\draw [dashed, ultra thick] (0,0) -- (0,3);
\draw [line width=0.10cm] (5,0) -- (5,3);
\draw [line width=0.10cm] (2,1) -- (2,3);
\draw [line width=0.10cm] (0,1) -- (2,1);
\draw [line width=0.10cm] (0,0) -- (5,0);
\draw [dashed, ultra thick] (0,3) -- (5,3);
\node at (1,2) {B'};
\node at (3.5,0.8) {B''};
\node at (2.5, -0.5) {Figure 4};
\end{tikzpicture}

\hskip .5in
\begin{tikzpicture} [scale=1]
\draw [dashed, ultra thick] (0,0) -- (0,3);
\draw [line width=0.10cm] (5,0) -- (5,3);
\draw [dashed, ultra thick] (3,1.5) -- (5,1.5);
\draw [dashed, ultra thick] (3,0) -- (3,1.5);
\draw [line width=0.10cm] (0,0) -- (5,0);
\draw [dashed, ultra thick] (0,3) -- (5,3);
\node at (1,2) {B''};
\node at (3.5, 0.8) {B'};
\node at (2.5, -0.5) {Figure 5};
\end{tikzpicture}}
\newcommand{\mynewpicture}[1][ ]{
\begin{tikzpicture} [scale=1]
\draw [dashed, ultra thick] (0,0) -- (0,3);
\draw [line width=0.10cm] (5,0) -- (5,3);
\draw [line width=0.10cm] (2,1) -- (2,3);
\draw [line width=0.10cm] (0,1) -- (2,1);
\draw [line width=0.10cm] (0,0) -- (5,0);
\draw [dashed, ultra thick] (0,3) -- (5,3);
\node at (1,2) {B'};
\node at (3.5,0.8) {B''};
\node at (2.5, -0.5) {Figure 4};
\end{tikzpicture}
\hskip .5in
\begin{tikzpicture} [scale=1]
\draw [dashed, ultra thick] (0,0) -- (0,3);
\draw [line width=0.10cm] (5,0) -- (5,3);
\draw [dashed, ultra thick] (3,1.5) -- (5,1.5);
\draw [dashed, ultra thick] (3,0) -- (3,1.5);
\draw [line width=0.10cm] (0,0) -- (5,0);
\draw [dashed, ultra thick] (0,3) -- (5,3);
\node at (1,2) {B''};
\node at (3.5, 0.8) {B'};
\node at (2.5, -0.5) {Figure 5};
\end{tikzpicture}}

\newcommand{\mypicture}[1][ ]{
\begin{tikzpicture} [scale=1]
\draw [line width=0.10cm] (0,0) -- (5,0);
\draw [line width=0.10cm] (5,0) -- (5,3);
\draw [dashed, ultra thick] (0,3) -- (5,3);
\draw [dashed, ultra thick] (0,0) -- (0,3);
\draw [line width=0.10cm] (3,0) -- (3,3);
\node at (1.5,1.5) {B1};
\node at (4,1.5) {B2};
\node at (2.5, -0.5){Figure 2};
\end{tikzpicture}
\hskip  .5in
\begin{tikzpicture}[ ] [scale=0.8]
\draw [line width=0.10cm] (7,0) -- (12,0);
\draw [line width=0.10cm] (12,0) -- (12,3);
\draw [dashed, ultra thick] (7,0) -- (7,3);
\draw [dashed, ultra thick] (7,3) -- (12,3);
\draw [line width=0.10cm] (7,1) -- (12,1);
\node at (9.5,2) {B1};
\node at (9.5,0.5) {B2};
\node at(9.5, -0.5) {Figure 3};
\end{tikzpicture}}

For  $a''<a' <a$,  $b <b' <b''$  and 
$B": =(a",a]\times [b,b"),$   $B': =(a',a]\times [b,b),$
the inclusion of vector spaces $ (\mathbb F^f_r(a'',b)+\mathbb F^f_r(a,b''))  \subseteq  (\mathbb F^f_r(a',b)+\mathbb F^f_r(a,b'))$  induces the canonical 
surjective linear map $\pi^{B'}_{B'', r}: \mathbb F^f_r(B'') \to \mathbb F^f_r(B').$  
\vskip .1in 

For  $0 <\epsilon' < \epsilon$  consider $B(a,b;\epsilon ')\subset (a-\epsilon, a]\times [b, b+\epsilon')= B_1 \subset B(a,b; \epsilon) $ and $B(a,b;\epsilon ')\subset (a-\epsilon', a]\times [b, b+\epsilon)= B_2 \subset B(a,b; \epsilon) .$ 
One has 

\hskip 1.5in $\pi^{\epsilon'}_{\epsilon,r} = \pi _{B_1,r}^{B(a,b; \epsilon')}\cdot \pi _{B(a,b;\epsilon),r}^{B_1} = \pi _{B_2,r}^{B(a,b; \epsilon')}\cdot \pi _{B(a,b;\epsilon),r}^{B_2} .$
\vskip .1in 
Consider the diagram 
\hskip .5in $\xymatrix { &\mathbb F^f_r(a,b)\ar [ld]_{\pi^{\epsilon}_{(a,b),r}} \ar[rd]^{\pi^{\epsilon'}_{(a,b),r}}
 &\\
\mathbb F^f_r(B(a,b;\epsilon))\ar[rr]^{\pi_\epsilon ^{\epsilon'}}&&\mathbb F^f_r(B(a,b;\epsilon')) }$

\noindent and denote  by $\hat\delta^f_r(a,b)$ and $\pi_r(a,b)$ the vector space $$\boxed{\hat \delta_r(a,b) = \varinjlim_{\epsilon\to 0} \mathbb F^f_r(B(a,b;\epsilon))}$$ 
and the surjective linear map 
$$\boxed{\pi_r (a,b): \mathbb F^f_r(a,b) \to \hat \delta _r(a,b) =  \varinjlim_{\epsilon\to 0} \pi_{(a,b)}^\epsilon}.$$ 
The space $\hat \delta_r(a,b) $   is  of finite dimension since $\mathbb F^f_r(a,b)$ is 
and denote this dimension by  $$ \boxed{\delta^f_r(a,b)= \dim \hat \delta^f_r(a,b)}.$$
\vskip .1in

In view of Observation \ref{O31} item 4  one proposes the following definition.  
 \begin{definition}\label {D1}\
 
A real number $t$  is  a {\bf homologically regular value} (w.r. to the field $\kappa$) if there exists $\epsilon (t) >0$  s.t. for any  $0<\epsilon <\epsilon(t)$
the inclusions $\mathbb I_{t-\epsilon}(r)\subseteq\mathbb I_t(r)\subseteq  \mathbb I_{t+\epsilon}(r)$  and 
$\mathbb I^{t-\epsilon}(r)\supseteq\mathbb I^t(r)\supseteq  \mathbb I^{t+\epsilon}(r)$ are equalities and {\bf homologically critical value} if not a {\bf homological regular value}
and let  $CR(f)$ be the set of  homological critical critical values. 
\end{definition}
By Observation \ref{O31} item 5,  $f$  weakly tame  
implies $CR(f)\subseteq Cr(f).$ 

\begin{obs} \label {P42}(cf \cite {B1})\ 

If $X$ is an ANR and $f$ is a continuous proper map then  
$CR(f)$ is a discrete set. 
 
If  $\delta^f_r(a,b)\ne 0$ then $a,b\in CR(f).$

If $f$ is tame and $\delta^f_r(a,b)\ne 0$ then $\hat \delta^f_r(a,b)= \mathbb F^f_r(B(a,b;\epsilon))$  for 
any $\epsilon <\epsilon (f).$\footnote {this observation holds also for $f$ continuous with the appropriate definition of $\epsilon (f).$}
\end{obs}
\vskip .1in

{\bf The configurations in the case of a real valued map}
\vskip .1in 

Suppose $f:X\to \mathbb R$ with $X$ compact  ANR and $f$ continuous. The assignment $\mathbb R^2\ni (a,b) \rightsquigarrow \delta^f_r(a,b)$ defined above is a {\it configuration of points} in $\mathbb R^2\equiv \mathbb C,$ which determines and is determined by  a monic polynomial $P^f_r(z)$  whose roots are the points in the support of $\delta^f_r$ with multiplicities the values of $\delta^f_r,$ and the assignment $\hat\delta^f_r$ is a {\it configuration of vector spaces} 
which  refines $\delta^f_r.$

If $\kappa= \mathbb R$ or $\mathbb C$ and $H_r(X)$ is equipped with a scalar product then the canonical splitting $s_r(a,b): \hat\delta^f_r(a,b)\to \mathbb F_r(a,b)$ of $\pi_r(a,b):\mathbb F_r(a,b) \to \hat\delta^f_r(a,b)$  given by the orthogonal complement of $\ker \pi_r(a,b)$ realizes $\hat\delta^f_r(a,b)$ as a subspace 
$$\boxed{\hat{\hat \delta }^f_r(a,b):= s_r(a,b)(\hat\delta ^f_r(a,b))\subseteq \mathbb F_r(a,b)\subseteq H_r(X)}.$$ 

It turns out that the points $(a,b)\in \supp \ \delta^f_r$ with $a\leq b$ are exactly the closed $r-$bar codes $[a,b]$  and  with $a>b$ are exactly the $(r-1)-$open bar codes $(b,a)$ defined in \cite{CSD09} and  \cite{BD11} for the level persistence of $f.$  
\vskip .1in

Note:  One can view the configurations $\delta^f_r$ 
and 
$\hat{\hat \delta} ^f_r$ 
 in analogy with 
 the configuration $\delta^T$ of eigenvalues with multiplicity, and the configuration $\hat{\hat \delta}^T$ of corresponding generalized eigenspaces, associated to 
a linear map $T:V\to V,$ $V$ a finite dimensional complex vector space.  
 The comparison provides  remarkable similarities which deserve to be inspected in case of a compact smooth Riemannian manifold and a Morse function. 
 
 \vskip .1in
{\bf The configurations in the case of an angle valued map}
\vskip .1in 

Suppose $f:X\to \mathbb S^1$ with $X$ compact.  Let $\tilde f: \tilde X\to \mathbb R$ be an infinite cyclic cover of $f,$ and  consider the homeomorphism $\tau : \tilde X\to \tilde X$ provided by the positive generator of the group of deck transformation $\mathbb Z;$  hence $\tilde f\cdot \tau = \tilde f + 2\pi.$ 
 The map $\tau$ induces the isomorphism $T_r: H_r(\tilde X)\to H_r(\tilde X)$  which restricts to  $T_r: \mathbb F_r^{\tilde f}(a,b)\to \mathbb F^{\tilde f}_r(a+2\pi, b+2\pi)$ and induces the isomorphism $T_r:\hat \delta^f_r(a,b)\to \hat\delta^f_r(a+2\pi, b+2\pi).$
\vskip .1in 
Consider the quotient space  $\mathbb T:= \mathbb R^2/\mathbb Z$  
identified to $\mathcal C\setminus 0$ 
by  
$\langle a, b\rangle \rightarrow e^{ia + (b-a)},$ cf subsection 2.1 and define  
$$\boxed{\delta^f_r (\langle a,b\rangle)= \delta^f_r(z):= \delta^{\tilde f}_r(a,b)},$$
 $$\boxed {\hat \delta ^f_r(\langle a,b\rangle)= \hat \delta ^f_r(z):= \oplus_{n\in \mathbb Z} \mathbb F^{\tilde f}_r(a+2n\pi, b+2n\pi)}$$  and 
 $$T_r(\langle a, b\rangle)= \oplus _{n\in \mathbb Z} T_r(a +2n\pi, b+2n\pi) : \hat\delta ^f_r(\langle a, b\rangle)\to \hat\delta ^f_r(\langle a, b\rangle).$$ The pair $({\hat \delta}^f_r(\langle a, b\rangle), T_r (\langle a, b\rangle)$ defines a  $\kappa[t^{-1}, t]-$module which is free 

It turns out that the points $e^{ia +(b-a)} \in \supp \delta^f_r$ with $a\leq b, a\in [0, 2\pi)$ are exactly the closed $r-$bar codes $[a,b]$ 
and  with $a>b, a\in [0,2\pi)$ are exactly the the $(r-1)-$open bar codes $(b,a)$ defined in \cite{BD11}.  
\vskip .1in 

As already pointed out in subsection 2.3, when $\kappa = \mathbb C$  the  algebra $\mathbb C[t^{-1},t]$ can be canonically completed to the finite von-Neumann algebra $L^\infty (\mathbb S^1).$ Additional data (for example a $\mathbb C[t^{-1},t]-$inner product on $H^N_r(X;\xi),$ or a Riemannian metric on $X$ when $X$ is a Riemannian manifold, or a triangulation of $X$ when $X$ is a simplicial complex) lead to a completion of $H^N_r(X;\xi)$  as Hilbert $L^\infty(\mathbb S^1)-$module, the $L_2-$homology $H^{L_2}(\tilde X),$  and of $\hat \delta^f_r(\langle a, b\rangle)$as a closed Hilbert submodule 
of $H^{L_2}(\tilde X).$ 
The procedure of such completions is described in \cite {B2} section 2  and called the {\it von-Neumann completion}.

The assignments $\delta^f_r,$ $\hat \delta^f_r,$ and $\hat{\hat\delta}^f_r$ are configurations of points with multiplicities, free $\kappa[t^{-1},t]-$modules and  $L^{\infty}(\mathbb S^1)-$Hilbert modules respectively.
\vskip .2in 

{\bf The Jordan cells for an angle valued map}
\vskip .1in 

For $f:X\to \mathbb S^1$ tame and $\theta\in \mathbb S^1$ denote by $X_\theta:=f^{-1}(\theta)$ and by $\overline X_\theta$ the two sided compactification of $f^{-1} (\mathbb S^1\setminus \theta)$ by $f^{-1}(\theta),$  called in \cite{B3} the {\it cut of $f$ at $\theta$}. The space $\overline X_\theta$ is homeomorphic to the compact space $\tilde f^{-1} [t, t+2\pi]$ for any $t\in \mathbb R$ with $p(t)=\theta.$ The inclusions 
$$\xymatrixcolsep{4pc}\xymatrix{X_\theta= \tilde f^{-1}(t)\ar[r]^{\subset}&\overline X_\theta = f^{-1}([t, t+2\pi] & X_\theta= f^{-1}(t+2\pi)\ar[l]_{\supset}}$$ induce in homology the linear map 
$$ \xymatrix{H_r(X_\theta)\ar[r]^{a}& H_r(\overline X_\theta)& H_r(X_\theta)\ar[l]_b}$$  which can be regarded as a {\it linear relation}, cf \cite {BH}, \cite {B3}, or as a graph representation of the oriented graph $G_2.$  

The oriented graph $G_2$ has two vertices $v , w$ and two oriented edges from $v$ to $w$ 
 denoted by $\alpha$ and $\beta$ as indicated below
 
 \hskip 2.5 in  \xymatrixcolsep{5pc} 
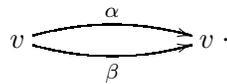
\xymatrix{v \ar@/^/[r]^\alpha \ar@/_/[r]_\beta & v}.
 
A linear representation $\rho$ of $G_2$ is provided by f.d. two vector spaces $V$ and $W$ associated to $v$ and $w$ and two linear maps $a,b: V\to W$ associated to the edges $\alpha, \beta$. The concept of isomorphism of representations direct sum of representations and indecomposable representations are obvious and, as in the case of an arbitrary finite oriented graph, each representation has a decomposition as  sum of a unique (up to isomorphism) collection of indecomposables; the decomposition is not unique.  
If $\kappa$  is algebraically closed the list of {\it indecomposables} can be recovered  from an old theorem of Kronecker (a proof of Kronecker theorem  can be found  in  \cite {Be})  and is provided below. 

\begin {enumerate}
\item 
Representation denoted by $\rho^+(r)$ has $V=\kappa^{r}, W=\kappa^{r+1},$
$a=\begin{bmatrix} 
Id_r\\0
\end{bmatrix},$
$b=\begin{bmatrix} 
0\\Id_r
\end{bmatrix}.$

\item 
Representation denoted by $\rho^-(r)$ has $V=\kappa^{r+1}, W=\kappa^{r},$
$a=\begin{bmatrix} 
Id_r&0
\end{bmatrix}, $
$b=\begin{bmatrix} 
0&Id_r
\end{bmatrix}.$

\item 
Representation denoted by  $(\lambda, k)$  called Jordan cells, has $V=\kappa^r, W=\kappa^{r},$
$a= T(\lambda,k),$ $b=\begin{bmatrix} 
Id_k\\0
\end{bmatrix}.$
 \end{enumerate}  

One defines the set  $\mathcal J_r(f,\theta)$ the collection of the Jordan cells associated to the $G_2$ representation given by 
 $$ \xymatrix{H_r(X_\theta)\ar[r]^{a}& H_r(\overline X_\theta)& H_r(X_\theta)\ar[l]_b}.$$

\section {The results}


As notices in Section 2  the configuration  $\delta^f_r$ defined in Section 3 can be equally regarded as a monic polynomial $P^f_r(z)$ whose zeros are the complex numbers $z\in \supp \ \delta^f_r$ with multiplicities equal to $\delta ^f_r(s).$  
\vskip .1in 
{\bf  Results about real valued maps}

\begin{theorem} (Topological results) \label {T1}\ 

Suppose $X$ compact and $f:X\to \mathbb R$ continuous. Then the following holds.
\begin{enumerate}
\item  
If $P^f_r(z)=0,$ equivalently $\delta^f_r(z)\ne 0$  with $z=(a+i b),$ then $a,b \in CR(f).$
\item The configuration $\delta^f_r\in \mathcal C_{\dim H_r(X)}(\mathbb C),$ the configuration $\hat \delta^f_r$ satisfies  $\oplus_{z\in \mathbb C} \hat \delta f_r(z)\simeq H_r(X)$ and if $\kappa= \mathbb R\  \rm{or}\  \ \mathbb C$ and $H_r(X)$ is equipped with a Hilbert space structure (i.e. a scalar product) then the configuration $\hat{\hat \delta} ^f(r)\in \mathcal C_{H_r(X)}(\mathbb C)$ and satisfies $\hat{\hat\delta}^f_r(z)\perp \hat{\hat\delta}^f_r(z)$ for $z\ne z'.$
\item For $f$ in an open and dense subset of the space of continuous real valued maps 
equipped with the compact open topology one has 
$\delta_r^f(z)= 0$ or $1.$
\end{enumerate}
\end{theorem}

\vskip .1in 
\begin{theorem} (Stability)\label{T2}\

Suppose $X$ is a 
compact ANR. 

1. The assignment $ f 
\rightsquigarrow \delta^f_r$ provides a continuous  map from the space of real valued maps 
equipped with the compact open topology to 
the space of configurations $\mathcal C_{b_r} (\mathbb R^2) = \mathcal C_{b_r}(\mathbb C)\simeq \mathbb C^{b_r},$\  $ b_r=\dim H_r(X),$ equipped with the collision topology, equivalently to the space of monic polynomials of degree $b_r.$ 

Moreover, with respect to the canonical metric $\underline D$  (cf Observation \ref {O211}) on the space of configurations $\mathcal C_{b_r}(\mathbb R^2),$ and the metric $D(f,g):= || f-g||_\infty= sup _{x\in X} |f(x)- g(x)|$ on the space of continuous maps  one has $$\underline D (\delta^f , \delta^g) < 2 D(f,g).$$  

2. If $\kappa= \mathbb R$ or $C,$ and $H_r(X)$ are equipped with scalar products  then the  assignment $f\rightsquigarrow \hat{\hat \delta}^f_r$ is also continuous provided that  $\mathcal C_{H_r(X)}(\mathbb C)$ is equipped with the collision topology described in subsection 2.1. 
\end{theorem}

\vskip .1in

\begin{theorem}  (Poincar\'e Duality) \label{T3}\

 Suppose $X$ is a closed topological manifold of dimension $n$ which is $\kappa-$orientable and $f: X\to \mathbb R$ a continuous map. Then the following holds.

1.  $\delta^f_r(a,b)= \delta^{f}_{n-r}(b,a).$ 

2. If $\kappa= \mathbb R, \mathbb C$ and the vector spaces $H_r(X)'$s  are equipped with scalar products  then the canonical isomorphism induced by the Poincar\'e duality and the scalar products,  
$PD_r: H_r(X)\to H_{n-r}(X),$  intertwines the configuration $\hat{\hat \delta}^f_r$ and $\hat{\hat\delta}^{f}_{n-r}\cdot \tau$ \ where $\tau (a,b)= (b,a).$ In particular if $X$ is a closed Riemannian manifold, hence $H_r(X)$ identifies to the space of harmonic $(n-r)-$ differential forms, then the Hodge star operator intertwines  $\hat{\hat \delta}^f_r$ with $\hat{\hat\delta}^{f}_{n-r}\cdot \tau.$ \end{theorem} 
\vskip .1in
 {\bf  Results about angle  valued maps}
\vskip .1in

Let $f:X\to \mathbb S^1$ be a continuous map, $X$ compact ANR,  and let $\xi:= \xi_f\in H^1(X;\mathbb Z)$ be the integral cohomology class represented by $f.$  Let $\tilde X\to X$ be an infinite cyclic cover associated to $\xi.$  
If $\kappa=\mathbb C$ let $H_r^{L_2}(\tilde X)$  be the von--Neumann completion of $H^N(X;\xi)$ 
as described in \cite{B2}.

\begin{theorem} (Topological results) \label {T4}\ 

Suppose $X$ compact ANR and $f:X\to \mathbb S^1$ continuous map. Then the following holds.

\begin{enumerate}
\item   
If $P^f_r(z)=0,$ equivalently  $\delta_r^f(z)
\ne 0$ with  $z= e^{ia +(b-a)},$ then $e^{ia}, e^{ib} \in CR(f)$   $(e^{ia}, e^{ib}\in \mathbb S^1).$
\item 
 The configuration  $ \delta^f_r(z)\in \mathcal C_{\beta^N_r(X;\xi_f)} (\mathbb C\setminus 0),$  the configuration $\hat \delta^f_r$ satisfies $\oplus \hat\delta^f_r\simeq H^N_r(X;\xi)$ and  if $\kappa =\mathbb C$ then the configuration 
$\hat{\hat \delta} ^f_r\in \mathcal C_{H^{L_2}(\tilde X)}(\mathbb C\setminus 0)$ and satisfies  $\hat{\hat \delta} ^f_r(z)\perp \hat{\hat \delta} ^f_r(z')$ for $z\ne z'.$
\item If $C_\xi(X, \mathbb S^1)$ denotes the set of continuous maps in the homotopy class determined by $\xi$  equipped with the compact open topology then for $f$ in an open and dense subset  of maps of $C_\xi(X,\mathbb S^1) $ one has $\delta^f(z)= 0\ \rm {or} \ 1.$
\end{enumerate}
\end{theorem}

\vskip .1in 
\begin{theorem} (Stability)\label{T5}

Suppose $X$ is a 
compact ANR and  $\xi\in H^1(X;\mathbb Z).$  Then the following holds.

1. The assignment 

\hskip .7in $ C(X, \mathbb S^1)_\xi \ni f 
\rightsquigarrow \delta^f_r\in \mathcal C_{\beta^N_r(X;\xi)}(\mathbb C\setminus 0)\equiv \mathcal C_{\beta^N_r(X;\xi)}(\mathbb R^2)$

equivalently 

\hskip .9in $C(X, \mathbb S^1)_\xi \ni f   
\rightsquigarrow P_r^f(z)\in \mathbb C^{\beta_r^N(X;\xi)}\times (\mathbb C\setminus 0)$

provides  a  continuous  map from $C_\xi(X, \mathbb S^1),$ the set of continuous maps in the homotopy class determined by $\xi$  equipped with the compact open topology, to the space 
of configurations $\mathcal C_{\beta^N_r(X;\xi)}(\mathbb C\setminus 0)$ equivalently $\mathbb C^{\beta_r^N(X;\xi)}\times (\mathbb C\setminus 0).$ 

Moreover, with respect to the canonical metric $\underline D$ on $C_{\beta^N_r(X;\xi)}(\mathbb T)$ and the complete metric $D$ on the space $C_\xi(X, \mathbb S^1)$ given by $D(f,g):=  sup _{x\in X} d(f(x),g(x)),$ $d$ the distance on $\mathbb S^1= \mathbb R/ 2\pi \mathbb Z,$ one has 
$$\underline D (\delta^f , \delta^g) < 2\pi D(f,g).$$  
2. If $\kappa= \mathbb C$ and the space of configurations $\mathcal C_{H_r^{L_2}(\tilde X)}(\mathbb C\setminus 0)$ is equipped with the  collision topology 
then the assignment $f \rightsquigarrow \hat{\hat \delta}^f_r$  is continuous. 
\end{theorem}

\begin{theorem}  (Poincar\'e Duality) \label{T6}

Suppose $M$ is a closed topological manifold of dimension $n$ which is $\kappa-$orientable and $f: M\to \mathbb S^1$  is a continuous map.  
 Then one has
\begin{enumerate}
\item  
$\delta^f_r(\langle a,b\rangle )= \delta^{ f}_{n-r}(\langle b, a\rangle ),$ equivalently $\delta^f_r(z)= \delta^f_{n-r}(\tau z)$ with $\tau(z)= z^{-1} e^{i\ln |z|}.$ Here $\langle a, b \rangle$ denotes the element of $\mathbb T$ represented by $(a,b)\in \mathbb R^2.$
\item If $\kappa=\mathbb C$ and $M$ is a closed Riemannian manifold then the canonical isomorphism  of $H^{L_2}_r(\tilde M)$ to $H^{L_2}_{n-r}(\tilde M)$ induced by the Riemannian metric (via $L_2$ harmonic forms and the Hodge star operator) intertwines the configuration $\hat{\hat \delta}^f_r$ and $\hat{\hat \delta}^f_{n-r}\cdot \tau $ when regarded as configurations on $\mathbb R^2/\mathbb Z =\mathbb T.$

\end{enumerate}
\end{theorem} 

In Section 3  for a weakly tame map $f:X\to \mathbb S^1$  and an angle $\theta\in \mathbb S^1$ we have defined the  collection of Jordan cells $\mathcal J_r(f,\theta),$ all computable by effective algorithms. They have  the following properties.

\begin{proposition}\label {P37}\
\begin{enumerate} 
\item If $f:X\to \mathbb S^1$ is a weakly tame map then the set $\mathcal J_r(f,\theta)$ is independent on $\theta,$ so  the notation $\mathcal J_r(f;\theta)$ can be abbreviated to $\mathcal J_r(f).$
\item If $f_1:X_1\to \mathbb S^1$ and $f_2:X_2\to \mathbb S^1$  are two weakly tame maps and 
$\omega:X_1\to X_2$ a homeomorphism  s.t. $f_2\cdot \omega$ and $f_1$ are homotopic then $\mathcal J_r(f_1)=  \mathcal J_r(f_2).$ 
\end{enumerate} 
\end{proposition}
This permits to define for any pair $(X,\xi)$ with $X$ a space homotopy equivalent to a compact ANR 
 and $\xi\in H^1(X; \mathbb Z)$ the invariant $\mathcal J_r(X,\xi)$ by   $\mathcal J_r(X,\xi):=\mathcal J_r(f)$ where $f:Y\to \mathbb S^1$ is a simplicial  map defined on the simplicial complex  $Y$ homotopy equivalent to $X$  by a homotopy equivalence $\omega:X\to Y$ 
 s.t. $f\cdot \omega$ represents $\xi.$ 
 In view of the discussion on the topology of compact Hilbert cube manifolds such pairs $(Y, \omega)$ exist.
 The invariant $\mathcal J_r(X;\xi)$ satisfies the following.
\begin{theorem}\label {T7}\
\begin{enumerate}
\item If $\omega:X_1\to X_2$ is a homotopy equivalence s.t.  $\omega^\ast (\xi_2)= \xi_1,$ \ \ 
$\xi_1\in H^1(X_1,\mathbb Z),$ \ $\xi_2\in H^1(X_2,\mathbb Z)$ and $X_1$ and $X_2$ have the homotopy type of a compact ANR then $\mathcal J(X_1,\xi_1)= \mathcal J_r(X_2, \xi_2).$
\item If $X$ is a compact ANR then $\mathcal J_r(X,\xi)$ are exactly
the Jordan cells of the monodromy $$T_r(X,\xi): V_r(X,\xi) \to V_r(X,\xi).$$
\end{enumerate}
\end{theorem} 
\vskip .1in  

Introduce the set 
$$\mathcal J_r(X;\xi)(u):= \{ (\lambda, k)\in \mathcal J_r(X:\xi) \mid \lambda=u\}$$
and for a finite set $S$ denote by $\sharp S$ the cardinality of $S.$

For any field $\kappa$ one has the following relation between Betti numbers , Novikov Betti numbers and Jordan cells.

\begin{theorem}\
$\beta_r(X)= \beta^N_r(X,\xi) + \sharp \mathcal J_r(X,\xi)(1) + \sharp \mathcal J_{r-1}(X,\xi)(1).$
\end{theorem}

\section {About the proof}

The proof of Theorems \ref{T1}, \ref {T2} , \ref {T3}  is contained partially in \cite {BH} and as stated in \cite {B1},  of Theorems \ref {T4}, \ref{T5} , \ref{T6} partially in \cite {BH} and as stated in \cite {B2}, and of Proposition \ref{P37} and Theorem \ref{T7} in \cite{BH} and \cite {B3}.

The  proofs are done first for nice spaces (homeomorphic to simplicial complexes) and tame maps and then extended to an arbitrary compact ANR and arbitrary continuous map  based on results on compact Hilbert cube manifolds as summarized in Theorem \ref{T54} below.

As far as the first step is concerned the following propositions of various level of complexity are   essential intermediate results whose proofs  are contained in \cite {B1}.

\begin{proposition}\label {P1}\

Let $a'<a<a'' ,$  $b< b''$ and $B_1,$ $B_2,$ and $B$  the boxes $B_1= (a',a]\times [b,b''),$ $B_2= (a,a'']\times [b,b'')$
and $B= (a',a'']\times [b,b'') $ (see Figure 2).

1. The inclusions $B_1\subset B$ and $B_2\subset B$ induce the linear maps 
$ i_{B_1,r} ^B: \mathbb F_r(B_1)\to  \mathbb F_r(B)$ and $\pi_{B,r}^{B_2}:  \mathbb F_r(B)\to  \mathbb F_r(B_2)$  such that the following sequence is exact
$$\xymatrix{0\ar [r]& \mathbb F_r(B_1)\ar[r]^{i_{B_1,r}^B} & \mathbb F_r(B)\ar[r]^{\pi_{B,r}^{B_2}} & \mathbb F_r(B_2)\ar[r]& 0}.$$ 

2. If $\kappa = \mathbb R$ or $\mathbb C$ and $H_r(X)$ is equipped with a scalar product hence $\mathbb F_r(B)'$s are canonically realized as subspaces 
${\bf H}_r(B)\subseteq H_r(M)$ then   
$${\bf H}_r(B_1) \perp {\bf H}_r(B_2)$$ and
$${\bf H}_r(B)= {\bf H}_r(B_1) + {\bf H}_r(B_2).$$
\end{proposition}

\begin{proposition}\label{P2}\

Let $a'<a,$ $b'<b< b''$ and $B_1,$ $B_2,$ and $B$  the boxes $B_1= (a',a]\times [b,b''),$ $B_2= (a,a'']\times [b',b)$
and $B= (a',a]\times [b',b'') $ (see Figure 3). 

1. The inclusions $B_1\subset B$ and $B_2\subset B$ induce the linear maps 
$i_{B_1,r}^B:  \mathbb F_r(B_1)\to  \mathbb F_r(B)$ and $\pi_{B,r}^{B_2}:  \mathbb F_r(B)\to  \mathbb F_r(B_2)$  such that the following sequence is exact
$$\xymatrix{0\ar [r]& \mathbb F_r(B_1)\ar[r]^{i_{B_1,r}^B} & \mathbb F_r(B)\ar[r]^{\pi_{B,r}^{B_2}} & \mathbb F_r(B_2)\ar[r]& 0}.$$ 

2. If $\kappa= \mathbb R$ or $\mathbb C$ and $H_r(X)$ is equipped with a scalar product then 
$${\bf H}_r(B_1) \perp {\bf H}_r(B_2)$$ and
$${\bf H}_r(B)= {\bf H}_r(B_1) + {\bf H}_r(B_2).$$ 
\end{proposition}

\mypicture {}

\vskip .1in
 \begin{proposition} \label {P08} (cf \cite{BH} Proposition 5.6)

Let $f:X\to \mathbb R$ be a  tame map and $\epsilon <\epsilon(f)/3.$ For any  map $g:X\to \mathbb R$ which satisfies  $|| f- g ||_\infty <\epsilon$ and $a,b\in Cr(f)$ critical values  one has
\begin{equation}\label{E3}
 \quad   \sum_{x\in D(a,b;2\epsilon)} \delta^g_r(x)=  \delta ^f_r(a,b), 
\end{equation}
\begin{equation} \label{E40}
 \quad \quad \supp \ \delta^{g}_r\subset \bigcup  _{(a,b)\in \supp\ \delta^{f}_r} D(a,b;2\epsilon). 
\end{equation}
If $\kappa= \mathbb R$ or $\mathbb C$ and in addition $H_r(X)$ is equipped with a scalar product the above statement can be strengthen to

\begin{equation}\label{E5}
x\in D(a,b;2\epsilon)\Rightarrow \hat{\hat \delta}^g_r(x)\subseteq \hat{\hat \delta} ^f_r(a,b), \ 
\quad   \oplus_{x\in D(a,b;2\epsilon)}\hat \delta^g_r(x)= \hat \delta ^f_r(a,b). 
\end{equation}
\end{proposition}
\vskip .1in

Theorems \ref{T1} and \ref {T3} follows essentially from the first two propositions which imply that 
$F_r$ is a measure on the sigma algebra generated by {\it boxes} with $\delta^f_r$ the {\it measure density}. Theorems \ref{T2} and \ref{T5}, in case the source of the map $f$ is a simplicial complex, uses essentially Proposition \ref{P08} and Theorems \ref{T3} and \ref{T6} use manipulation of Poincar\'e duality and alternative definition of $\hat \delta^f_r$ , cf \cite{B1}.  In case of Theorem \ref{T6} a more elaborated manipulation involving Poincar\'e duality for the open manifold $\tilde M,$ the infinite cyclic cover of $f: M\to \mathbb S^1,$ and the description of the torsion of the $\kappa[t^{-1},t]$ module $H_r(\tilde M)$  are  needed, cf \cite {B2}.
The proof of Theorem \ref{T7} involves the recognition of what in \cite{BH} and \cite{B3} is referred to as the {\it regular part of the linear relation} defined by the pair of linear maps $$ \xymatrix{H_r(X_\theta)\ar[r]^{a}& H_r(\overline X_\theta)& H_r(X_\theta)\ar[l]_b}.$$

Concerning the results about compact Hilbert cube manifolds used in this work  
recall that: 

The Hilbert cube $Q$ is the infinite product $Q =  \prod_{i\in \mathbb Z_{\geq 0}} I_i$ with $I_i= [0,1]; $ its  topology  is also  given by the metric $d(\overline u ,\overline v)= \sum _i |u_i- v_i|/ 2^i$ with $\overline u= \{u_i\in I_i, i\in \mathbb Z_{\geq 0}\}$  and  $\overline v= \{v_i\in I_i, i\in \mathbb Z_{\geq 0}\}.$ The space $Q$ is a compact ANR and so is any $X\times Q$ for $X$  any compact ANR. 

A compact Hilbert cube manifold is a compact Hausdorff space locally  homeomorphic to the Hilbert cube and is a compact ANR.
The following basic  results about Hilbert cube manifolds  can be found in \cite {CH}.

\begin{theorem} \label {T54}\ 
\begin{enumerate} 
\item (R Edwards) $X$ is a compact ANR iff $X\times Q$ is a compact Hilbert cube manifold.
\item (T.Chapman) Any  compact Hilbert cube manifolds is homeomorphic to $K\times Q$ \ for some finite simplicial complex $K.$ 
\item (T Chapman) If $\omega:X\to Y$ is a simple homotopy equivalence between two finite simplicial complexes with Whitehead torsion $\tau(\omega)=0$ then there exists a homeomorphism $\omega': X\times Q\to Y\times Q$ s.t. 
$\omega'$ and $\omega\times id_Q$ are homotopic. If $\omega$ is only a homotopy equivalence the same conclusion holds for $Q$ replaced by $Q\times \mathbb S^1.$ \footnote{some partial but relevant results on the line of Theorem \ref{T54} were due to J West as indicated in \cite {CH}}  
\end{enumerate}
\end{theorem}

If one writes $I^\infty= I^k\times I^{\infty-k}$ observe that 
 given $\epsilon >0$  for any  continuous  real or angle valued map $f$ defined on $K\times Q,$ $K-$ a simplicial complex, there exists $N$ large enough such that  $f$ is $\epsilon-$closed to
$g\cdot \pi,$  $\pi: K\times I^\infty\to K\times I^N$ the canonical projection  with $g$ a simplicial map defined on $K\times I^N.$ In particular any compact Hilbert cube manifold is a good ANR.
It can be also verified  by using the definitions that if $f:X\to \mathbb R\ \rm{or}\  \mathbb S^1$ is a continuous map $K,$ a compact ANR  and $f^K= f\times \pi,$ $\pi: X\times K\to X,$ then 
$\hat \delta^{f^K}(\langle a,b\rangle )= \oplus_{k\geq 0} \hat \delta^f_{r-k}(\langle a,b\rangle) \otimes H_k(K).$ 

\section {Some applications}

1. {\bf Geometric analysis}

Theorem \ref{T1} insures that a generic continuous  function provides one dimensional subspaces in homology with coefficients in a fixed field;  in particular for $\kappa=m\mathbb R$ or $\mathbb C$ for a closed Riemannian manifold a  generic continuous function provides an orthonormal base  (up to sign) 
in the space of harmonic forms.  

We  expect (but have not found this result in literature) that the eigenforms of the Laplace Beltrami operators for a generic Riemannian metric in any dimension provides a similar decomposition for the smooth differential forms  orthogonal to the harmonic forms. This is indeed the case in view of a result of Uhlenbeck\footnote{which claims that for a  closed  manifold equipped with a generic Riemannian metric the eigenvalues of the Laplace operator  are simple} for degree zero forms and for $n=2$ for all degree.  This shows that  a generic pair, Riemannian metric and smooth function,  provides an orthonormal  base up to sign (In Fourier sense) in the space of all differential forms;  in the same way  trigonometric functions on $\mathbb S^1$ provide an orthonormal  base (In Fourier sense) for smooth functions. This can be a useful tool in geometric analysis.
\vskip .1in 
2. {\bf Topology }
\begin {obs} \label {O53}\ 

1.  Theorem (\ref{T3})  implies that for a closed orientable manifold of dimension $n$ $(c,c')\in \supp \delta^f_r$ iff $(c',c)\in \supp \delta^f_{n-r}$ and both pairs appear with equal multiplicity $\delta^f_r(c,c')= \delta^f_{n-r}(c',c).$

2. Theorem \ref {T6} remains valid with  the same proof in case $M$ is a compact manifold  with boundary $(M,\partial M),$ 
provided $H^N_r(\partial M; \xi_{f_{\partial M}})$ \footnote {with $f_{\partial M}$ notation for the restriction of $f$ to $\partial M$}
vanishes for all $r.$ In particular, under the above vanishing hypothesis,  
$H^N_r(M;\xi_f)\simeq H^N_{n-r} (M;\xi_f).$
\end{obs}

\begin {corollary}\label {C}\

Suppose $(M^{2n},\partial M^{2n})$ is a compact orientable  manifold with boundary which has the homotopy type of a simplicial complex of dimension $\leq n$ and $\xi\in H^1(M;\mathbb Z)$ s.t $H^N_r(\partial M; \xi_{\partial M})=0$ for all $r.$
Then for any field $\kappa$ :

1.  $\beta^N_r(X:\xi)= \begin{cases}  0 \ \rm{if}\  r\ne n\\ (-1)^n\chi (M_n)\ \rm{if}\  r= n\end{cases},$ with $\chi (M)$ the Euler -Poincar\'e characteristic with coefficients in $\kappa.$

2. $\beta_r(X) =\begin{cases}  \alpha_{r-1} + \alpha_r \ \rm{if}\  r\ne n\\ \alpha_{n-1} + \alpha_n  +(-1)^n\chi (M_n)\ \rm{if}\  r= n\end{cases}, $
where $\alpha _r$ denotes the number of Jordan cells $(\lambda, k) \in J_r(M,\xi_f),$ with $\lambda=1.$  
\vskip.1in 
3. If $V^{2n-1}\subset M^{2n}$ is a  compact  proper sub manifold (i,e, $V\pitchfork \partial M, \footnote {$\pitchfork$= transversal}$ and $V\cap \partial M= \partial V$) representing a homology class in $H_{n-1}(M,\partial M)$  Poincar\'e dual to $\xi_f$  and $H_r(V)=0$ then the set of Jordan cells $J_r(M,\xi)$  is empty. 
\end{corollary}

Item 1. follows from Observation (\ref {O53}) and the fact that both Betti numbers and Novikov--Betti numbers calculate the same Euler--Poincar\'e characteristic.  Item 2  follows from Theorem 11 item c.  in \cite {BH}, and Item 3.  from the description of Jordan cells in terms   of linear relations  as provided in \cite{BH} or \cite{B3}.

As pointed out to us by L Maxim, the complement $X = \mathbb C^n\setminus V$ of a complex hyper surface $V\subset \mathbb C^n, V:= \{(z_1, z_2, \cdots z_n) \mid f(z_1, z_2, \cdots z_n)=0\}$ regular at infinity, 
equipped with the canonical class $\xi_f\in H^1(X; \mathbb Z)$ defined by $f: X\to \mathbb C\setminus 0$ is an open manifold with an integral  cohomology class  $\xi\in H^1(X;\mathbb Z)$  represented by $f/|f|: X\to \mathbb S^1.$ This manifold has compactifications to manifolds with boundary with  cohomology class which satisfies the hypotheses above.  
Item 1. recovers    
 a calculation of L Maxim, cf \cite{M14} and \cite {FM16} \footnote {The  Friedl-Maxim results  state the vanishing of more general and more sophisticated $L_2-$homologies and Novikov type  homologies. They can also be recovered via the appropriate Poincar\'e Duality type isomorphisms on similar lines.}   
that the complement of an algebraic hyper surface regular at infinity   
has vanishing Novikov homologies in all dimension but $n.$

\end{document}